\documentclass[12pt]{article}
\usepackage{amsmath}
\usepackage{latexsym}
\usepackage{epsfig}

\begin{document}

\begin{center}
\bigskip

\bigskip

{\Large Stochastic arbitrage return and its implication for option pricing }

\bigskip

{\large Sergei Fedotov and Stephanos Panayides}

\bigskip Department of Mathematics, UMIST, M60 1QD UK

\bigskip

Submitted to Physica A on 29 March 2004



\bigskip

\bigskip

\bigskip

\textbf{Abstract} \bigskip
\end{center}

The purpose of this work is to explore the role that random
arbitrage opportunities play in pricing financial derivatives. We
use a non-equilibrium model to set up a stochastic portfolio, and
for the random arbitrage return, we choose a stationary ergodic
random process rapidly varying in time. We exploit the fact that
option price and random arbitrage returns change on different time
scales which allows us to develop an asymptotic pricing theory
involving the central limit theorem for random processes. We
restrict ourselves to finding pricing bands for options rather
than exact prices. The resulting pricing bands are shown to be
independent of the detailed statistical characteristics of the
arbitrage return. We find that the volatility ''smile'' can also
be explained in terms of random arbitrage opportunities.
\pagebreak

\section{Introduction}

It is well-known that the classical Black-Scholes formula is consistent with
quoted options prices if different volatilities are used for different
option strikes and maturities \cite{Hull}. To explain this phenomenon,
referred to as the volatility ''smile'', a variety of models has been
proposed in the financial literature. These includes, amongst others,
stochastic volatility models \cite{Sircar,Lewis,Sircar2}, and the Merton
jump-diffusion model \cite{M}. Each of these is based on the assumption of
the absence of arbitrage. However, it is well known that arbitrage
opportunities always exist in the real world (see \cite{GA,SO}). Of course,
arbitragers ensure that the prices of securities do not get out of line with
their equilibrium values, and therefore virtual arbitrage is always
short-lived. One of the purposes of this paper is to explain the volatility
''smile'' phenomenon in terms of random arbitrage opportunities.

The first attempt to take into account virtual arbitrage in option pricing
was made by physicists in \cite{AE,KIL,KIS}. The authors assume that
arbitrage returns exist, appearing and disappearing over a short time scale.
In particular, the return from the Black-Scholes portfolio, $\Pi
=V-S\partial V/\partial S,$ where $V$ is the option price written on an
asset $S,$ is not equal to the constant risk-free interest rate $r$. In \cite
{KIL,KIS} the authors suggest the equation $d\Pi /dt=(r$ $+$ $x\left(
t\right) )\Pi ,$ where $x\left( t\right) $ is the random arbitrage return
that follows an Ornstein-Uhlenbeck process. In \cite{MO,MO2} this idea is
reformulated in terms of option pricing with stochastic interest rate. The
main problem with this approach is that the random interest rate is not a
tradable security, and therefore the classical hedging can not be applied.
This difficulty leads to the appearance of an unknown parameter, the market
price of risk, in the equation for derivative price \cite{MO,MO2}. Since
this parameter is not available directly from financial data, one has to
make further assumptions on it. An alternative approach for option pricing
in an incomplete market is based on risk minimization procedures (see, for
example, \cite{AS,BS,FM,Sch,Sc}).

In this paper we follow an approach suggested in \cite{Sircar2}
where option pricing with stochastic volatility is considered.
Instead of finding the exact equation for option price, we focus
on the pricing bands for options that account for random arbitrage
opportunities. We exploit the fact that option price and random
arbitrage return change on different time scales allowing us to
develop an asymptotic pricing theory by using the central limit
theorem for random processes \cite{FW}. The approach yields
pricing bands that are independent of the detailed statistical
characteristics of the random arbitrage return.

\section{Model with random arbitrage return}

Consider a model of $(S,B,V)$ market that consists of the stock, $S$, the
bond, $B$ and the European option on the stock, $V$. To take into account
random arbitrage opportunities, we assume that this market is affected by
two sources of uncertainty. The first source is the random fluctuations of
the return from the stock, described by the conventional stochastic
differential equation,
\begin{equation}
\frac{dS}{S}=\mu dt+\sigma dW,  \label{stock}
\end{equation}
where $W$ is a standard Wiener process. The second source of uncertainty is
a random arbitrage return from the bond described by
\begin{equation}
\frac{dB}{B}=rdt+\xi (t)dt,  \label{bond}
\end{equation}
where $r$ is the risk-free interest rate. Given that there are random
arbitrage opportunities, we introduce here the random process, $\xi (t),$
that describes the fluctuations of return around $rdt.$ The same mapping to
a model with stochastic interest rate is used in \cite{KIL,KIS,MO}. The
characteristic of the present model is that we do not assume that $\xi (t)$
obeys a given stochastic differential equation. For example, in \cite
{KIL,KIS,MO} $\xi (t)$ follows the Ornstein-Uhlenbeck process \cite{OK}.

It is reasonable to assume that random variations of arbitrage return $\xi
(t),$ are on the scale of hours. Let us denote this characteristic time by $%
\tau _{arb}.$ This time can be regarded as an intermediate one between the
time scale of random stock return (infinitely fast Brownian motion
fluctuations), and the lifetime of the derivative $T$ (several months): $0<<$
$\tau _{arb}<<T$. In what follows, we exploit this difference in time scales
and develop an asymptotic pricing theory involving the central limit theorem
for random processes.

Now we are in a position to derive the equation for $V.$ Let us consider the
investor establishing a zero initial investment position by creating a
portfolio $\Pi $ consisting of a long position of one bond, $B,$ $\frac{%
\partial V}{\partial S}$ shares of the stock, $S,$ and a short position in
one European option, $V,$ with an exercise price, $K,$ and a maturity, $T.$
The value of this portfolio is
\begin{equation}
\Pi =B+\frac{\partial V}{\partial S}S-V.  \label{portfolio}
\end{equation}

The Black-Scholes dynamic of this portfolio is given by two equations, $%
\partial \Pi /\partial t=0,$ and $\ \Pi =0$ (see \cite{Hull}). The
application of Ito's formula to (\ref{portfolio}) together with (\ref{stock}%
) and (\ref{bond}) with $\xi (t)=0,$ leads to the classical Black-Scholes
equation:
\begin{equation}
\frac{\partial V}{\partial t}+\frac{\sigma ^{2}S^{2}}{2}\frac{\partial ^{2}V%
}{\partial S^{2}}+rS\frac{\partial V}{\partial S}-rV=0.  \label{BS}
\end{equation}
The natural generalization of $\partial \Pi /\partial t=0$ is a simple
non-equilibrium equation:
\begin{equation}
\frac{\partial \Pi }{\partial t}=-\frac{\Pi }{\tau _{arb}},  \label{p}
\end{equation}
where $\tau _{arb}$ is the characteristic time during which the arbitrage
opportunity ceases to exist (see \cite{AE}).

By using a self-financing condition, $d\Pi =dV-\frac{\partial V}{\partial S}%
dS+dB,$ and Ito's lemma, one can derive from (\ref{stock}) and (\ref{bond})
the equation involving the option $V(t,S)$ and the portfolio value $\Pi (t,S)
$,
\begin{equation}
\frac{\partial V}{\partial t}+\frac{\sigma ^{2}S^{2}}{2}\frac{\partial ^{2}V%
}{\partial S^{2}}+rS\frac{\partial V}{\partial S}-rV+r\Pi +\xi (t)\Pi +\xi
(t)\left( S\frac{\partial V}{\partial S}-V\right) +\frac{\Pi }{\tau _{arb}}%
=0.  \label{basic}
\end{equation}
Note that this equation reduces to (\ref{BS}) when $\Pi =0$ and $\xi (t)=0.$

To deal with the forward problem we introduce the non-dimensional time
\begin{equation}
\tau =\frac{T-t}{T},\;\;\;0\leq \tau \leq 1,  \label{time}
\end{equation}
and a small parameter
\begin{equation}
\varepsilon =\frac{\tau _{arb}}{T}<<1.
\end{equation}
This parameter plays a very important role in what follows. In the limit $%
\varepsilon \rightarrow 0,$ the stochastic arbitrage return $\xi $ becomes a
function that is rapidly varying in time. It is convenient to use the
following notation, $\xi \left( \frac{\tau }{\varepsilon }\right) $ (see
\cite{Sircar2}). It follows from (\ref{p}) that the value of the portfolio $%
\Pi ,$ decreases to zero like $\exp (-t/\varepsilon T).$ Thus, one can
assume that $\Pi $ $=0$ in the limit $\varepsilon \rightarrow 0.$ We find
from (\ref{basic}) that the associated option price, $V^{\varepsilon }\left(
\tau ,S\right) ,$ obeys the following stochastic PDE
\begin{equation}
\frac{\partial V^{\varepsilon }}{\partial \tau }=\frac{\sigma ^{2}S^{2}}{2}%
\frac{\partial ^{2}V^{\varepsilon }}{\partial S^{2}}+rS\frac{\partial
V^{\varepsilon }}{\partial S}-rV^{\varepsilon }+\xi (\frac{\tau }{%
\varepsilon })\left( S\frac{\partial V^{\varepsilon }}{\partial S}%
-V^{\varepsilon }\right)   \label{basic2}
\end{equation}
subject to the initial condition:
\begin{equation}
V^{\varepsilon }(0,S)=\max (S-K,0)
\end{equation}
for a call option, where $K$ is the strike price. Here, for simplicity, we
keep the same notations for the non-dimensional volatility $\sigma $ and the
interest rate $r.$ The same PDE is used in \cite{KIL,KIS}, but with the
Ornstein-Uhlenbeck process for $\xi .$

\section{Asymptotic analysis: pricing bands for the options}

To analyze the stochastic PDE (\ref{basic2}), we have to specify the
statistical properties of the random arbitrage return $\xi $. Suppose that $%
\xi (\tau )$ is a stationary ergodic random process with zero mean, \ $<$ $%
\xi (\tau )$ $>=0,$ such that,
\begin{equation}
D=\int_{0}^{\infty }<\xi (\tau +s)\xi (\tau )>ds,  \label{dif}
\end{equation}
is finite. The key feature of this paper is that we do not assume an
explicit equation for $\xi (\tau )$ unlike the works \cite{KIL,KIS,MO},
where the random arbitrage return $\xi $ follows the Ornstein-Uhlenbeck
process.

According to the law of large numbers, $V^{\varepsilon }\left( \tau
,S\right) $ converges in probability to the Black-Scholes price, $%
V_{BS}\left( \tau ,S\right) ,$ as $\varepsilon \rightarrow 0.$ One can split
$V^{\varepsilon }(\tau ,S)$ into the sum of the Black-Scholes price, $%
V_{BS}\left( \tau ,S\right) $, and the random field $Z^{\varepsilon }(\tau
,S),$ with the scaling factor $\sqrt{\varepsilon },$ giving
\begin{equation}
V^{\varepsilon }\left( \tau ,S\right) =V_{BS}\left( \tau ,S\right) +\sqrt{%
\varepsilon }Z^{\varepsilon }\left( \tau ,S\right) .  \label{split2}
\end{equation}
Substituting (\ref{split2}) into (\ref{basic2}), and using the equation for $%
V_{BS}\left( \tau ,S\right) $, we get the following stochastic PDE for $%
Z^{\varepsilon }(\tau ,S)$,
\begin{equation}
\frac{\partial Z^{\varepsilon }}{\partial \tau }=\frac{1}{2}\sigma ^{2}S^{2}%
\frac{\partial ^{2}Z^{\varepsilon }}{\partial S^{2}}+(r+\xi (\frac{\tau }{%
\varepsilon }))(S\frac{\partial Z^{\varepsilon }}{\partial S}-Z^{\varepsilon
})+\frac{\xi (\frac{\tau }{\varepsilon })}{\sqrt{\varepsilon }}\left( S\frac{%
\partial V_{BS}}{\partial S}-V_{BS}\right) .  \label{error}
\end{equation}
Our objective here is to find the asymptotic equation for
$Z^{\varepsilon }(\tau ,S)$ as $\varepsilon \rightarrow 0.$ One
can see that Eq. ( \ref{error}) involves two stochastic terms
proportional to $\xi (\frac{\tau }{\varepsilon })$ and
$\varepsilon ^{-1/2}\xi (\frac{\tau }{\varepsilon })$. Ergodic
theory implies that the first term in its integral form converges
to zero as $\varepsilon \rightarrow 0,$ while the second term
converges weakly to a white Gaussian noise $\eta \left( \tau
\right) $ with the correlation function
\begin{equation}
<\eta (\tau _{1})\eta (\tau _{2})>=2D\delta \left( \tau _{1}-\tau
_{2}\right) ,  \label{delta}
\end{equation}
where the intensity of white noise, $D,$ is determined by
(\ref{dif}) (see \cite{FW} and Appendix A). Thus, in the limit
$\varepsilon \rightarrow 0,$ the random field $Z^{\varepsilon
}(\tau ,S),$ converges weakly to the field $Z(\tau ,S)$ that obeys
the asymptotic stochastic PDE:
\begin{equation}
\frac{\partial Z}{\partial \tau }=\frac{1}{2}\sigma ^{2}S^{2}\frac{\partial
^{2}Z}{\partial S^{2}}+r(S\frac{\partial Z}{\partial S}-Z)+\left( S\frac{%
\partial V_{BS}}{\partial S}-V_{BS}\right) \eta (\tau ),
\label{limitequation}
\end{equation}
with the initial condition $Z(0,S)=0$. This equation can be solved
in terms of the classical Black-Scholes Green function,
$G(S,S_{1},\tau ,\tau _{1}),$ to give
\begin{equation}
Z(\tau ,S)=\int_{0}^{\tau }\int_{0}^{\infty }G(S,S_{1},\tau ,\tau
_{1})\left( S_{1}\frac{\partial V_{BS}}{\partial S_{1}}-V_{BS}\left( \tau
_{1},S_{1}\right) \right) \eta (\tau _{1})dS_{1}d\tau _{1},  \label{solution}
\end{equation}
where
\begin{equation}
G(S,S_{1},\tau ,\tau _{1})=\frac{e^{-r(\tau -\tau _{1})}}{S_{1}\sqrt{2\pi
\sigma ^{2}(\tau -\tau _{1})}}e^{-\frac{[\ln (S/S_{1})+(r-\frac{\sigma ^{2}}{%
2})(\tau -\tau _{1})]^{2}}{2\sigma ^{2}(\tau -\tau _{1})}}
\label{Green}
\end{equation}
(see, for example, \cite{C}).

It follows from (\ref{solution}) that since $\eta \left( t\right) $ is the
Gaussian noise with zero mean, $Z\left( \tau ,S\right) $ is also the
Gaussian field with zero mean. The covariance
\begin{equation}
R\left( \tau ,S,Y\right) =<Z\left( \tau ,S\right) Z\left( \tau ,Y\right) >
\end{equation}
satisfies the deterministic PDE:
\begin{equation*}
\frac{\partial R}{\partial \tau }=\frac{1}{2}\sigma ^{2}S^{2}\frac{\partial
^{2}R}{\partial S^{2}}+\frac{1}{2}\sigma ^{2}Y^{2}\frac{\partial ^{2}R}{%
\partial Y^{2}}+r(S\frac{\partial R}{\partial S}+Y\frac{\partial R}{\partial
Y}-2R)+
\end{equation*}
\begin{equation}
2D\left( S\frac{\partial V_{BS}}{\partial S}-V_{BS}\left( \tau ,S\right)
\right) \left( Y\frac{\partial V_{BS}}{\partial Y}-V_{BS}\left( \tau
,Y\right) \right)   \label{PDE}
\end{equation}
with $R(0,S,Y)\equiv 0$ (see Appendix B). The pricing bands for the options
for the case of arbitrage opportunities can be given by
\begin{equation}
V_{BS}\left( \tau ,S\right) \pm 2\sqrt{\varepsilon U\left( \tau ,S\right) },
\label{band}
\end{equation}
where
\begin{equation}
U\left( \tau ,S\right) =R\left( \tau ,S,S\right) .  \label{UU}
\end{equation}
The variance $U\left( \tau ,S\right) $ quantifies the fluctuations around
the classical Black-Scholes price. It should be noted that it is independent
of the detailed statistical characteristics of the arbitrage return. The
only parameter we need to estimate is the intensity of noise $D$ which is
the integral characteristic of random arbitrage return (see (\ref{dif})).
One can conclude that the investor who employs the arbitrage opportunities
band hedging sells the option for
\begin{equation}
V_{BS}\left( \tau ,S\right) +2\sqrt{\varepsilon U\left( \tau ,S\right) },
\label{effective}
\end{equation}
where $U\left( \tau ,S\right) $ can be found from (\ref
{solution}) or (\ref{PDE}) to be
\begin{equation}
U\left( \tau ,S\right) =2D\int_{0}^{\tau }[\int_{0}^{\infty }G(S,S_{1},\tau
,\tau _{1})\left( S_{1}\frac{\partial V_{BS}}{\partial S_{1}}-V_{BS}\right)
dS_{1}]^{2}d\tau _{1}  \label{U}
\end{equation}
(see Appendix C).

One can see from (\ref{limitequation}) or (\ref{U}) that the large
fluctuations of $Z(\tau ,S)$ occur when the function $\left( S\frac{\partial
V_{BS}}{\partial S}-V_{BS}\right) $ takes the maximum value. That is, $S%
\frac{\partial ^{2}V_{BS}}{\partial S^{2}}=0$ (the Greek $\Gamma =0)$.

\section{Numerical results}

To determine how the random arbitrage opportunities affect the
option price, we solve equation (\ref{PDE}) numerically. Figure 1
shows a graph of the variance, $U(\tau ,S)=R\left( \tau
,S,S\right) $, as a function of both time, $\tau $, and asset
price, $S$. From this graph, we observe that the uncertainty
regarding the option value is greater for deep-in-the money
options. This finding is consistent with the empirical work given
in \cite {BCC}. In Figure 2, we plot the effective option price
given by (\ref{effective}) for $\varepsilon =0.1$, and compare it
with the Black-Scholes price. Note that deep-in-the money options
deviate the most from the Black-Scholes option price. As we move
near at-the-money options the deviation decreases.

Now we are in a position to discuss the ''smile'' effect, that is, the
implied volatility is not a constant, but varies with strike price $K$ and
time $\tau .$ Let us denote the implied Black-Scholes volatility by $%
\sigma ^{\varepsilon }\left( K,\tau \right).$ The formula
(\ref{effective}) for the effective option price implies
\begin{equation}
V_{BS}\left( \tau ,S;\sigma ^{\varepsilon }\left( K,\tau \right) ,K\right)
=V_{BS}\left( \tau ,S;\sigma ,K\right) +2\sqrt{\varepsilon U\left( \tau
,S;K\right) }.  \label{implied}
\end{equation}
This equation can be solved for $\sigma ^{\varepsilon }\left( K,\tau \right)
$ with $S$ and $\varepsilon $ fixed. It follows from (\ref{effective}) that $%
\sigma ^{\varepsilon }\left( K,\tau \right) \rightarrow \sigma $ as $%
\varepsilon \rightarrow 0.$ Fig.3 illustrates the `smile' effect,
when implied volatility $\sigma ^{\varepsilon }\left( K,\tau
\right) $ increases for deep-in and out-of-the money options.
Similar results are given in \cite
{Otto} where the characteristic time $\tau _{arb}$ depends on moneyness $%
\frac{S}{K}$. Numerical results give a similar smile curve where the implied
volatility becomes greater as we move away from at-the-money.

\section{Conclusions}

In this paper, we investigated the implications of random arbitrage return
for option pricing. We extended previous works by using a stationary ergodic
process for modelling random arbitrage return. We considered the case where
arbitrage return fluctuates on a different time scale to that of the option
price. This allowed us to use asymptotic analysis to find option pricing
bands rather than the exact equation for option value. We derived the
asymptotic equation for the random field that quantifies the fluctuations
around the classical Black-Scholes price and showed that it is independent
of the detailed statistical characteristics of the arbitrage return. In
particular, we showed that the risk from the random arbitrage returns is
greater for deep-in and out-of-the money options. This gives an explanation
of the `smile' effect observed in the market in terms of the random
arbitrage return.

\renewcommand{\theequation}{A-\arabic{equation}} \setcounter{equation}{0}

\section*{Appendix A}

The random process $%
x^{\varepsilon }(\tau )$ defined as
\begin{equation}
x^{\varepsilon }(\tau )=\frac{1}{\varepsilon ^{\frac{1}{2}}}\int_{0}^{\tau
}\xi (\frac{s}{\varepsilon })ds=\varepsilon ^{\frac{1}{2}}\int_{0}^{\frac{%
\tau }{\varepsilon }}\xi (s)ds  \label{new}
\end{equation}
converges weakly to the Brownian motion $B\left( \tau \right) $ as $%
\varepsilon \rightarrow 0$. It means that $\varepsilon ^{-1/2}\xi (\frac{%
\tau }{\varepsilon })$ converges weakly to a white Gaussian noise
$\eta \left( \tau \right) $ with the correlation function $<\eta
(\tau _{1})\eta (\tau _{2})>=$ $2D\delta \left( \tau _{1}-\tau
_{2}\right) $ (note that $\eta\left( \tau \right)  =\frac{d
B\left( \tau \right)
 }{d\tau }).$

The purpose of this Appendix is to show that
\begin{equation}
\lim_{\varepsilon \rightarrow 0}<[x^{\varepsilon }(\tau )]^{2}>=2D\tau .
\end{equation}
It follows from (\ref{new}) that
\begin{equation}
<[x^{\varepsilon }(\tau )]^{2}>=\varepsilon \int_{0}^{\frac{\tau }{%
\varepsilon }}\int_{0}^{\frac{\tau }{\varepsilon }}<\xi (s_{1})\xi
(s_{2})>ds_{1}ds_{2}.
\end{equation}
By using the well-known formula for stationary process $\xi (\tau )$ with
zero mean
\begin{equation}
\int_{0}^{\tau }\int_{0}^{\tau }<\xi (s_{1})\xi (s_{2})>ds_{1}ds_{2}=2\tau
\int_{0}^{\tau }<\xi (\tau +s)\xi (\tau )>ds-2\int_{0}^{\tau }s<\xi (\tau
+s)\xi (\tau )>ds,
\end{equation}
one finds that
\begin{equation}
E\{[x^{\varepsilon }(\tau )]^{2}\}=2\tau \int_{0}^{\frac{\tau }{\varepsilon }%
}<\xi (\tau +s)\xi (\tau )>ds-2\varepsilon \int_{0}^{\frac{\tau }{%
\varepsilon }}s<\xi (\tau +s)\xi (\tau )>ds.
\end{equation}
In the limit $\varepsilon \rightarrow 0$ the second term tends to zero, and
so we must have
\begin{equation*}
\lim_{\varepsilon \rightarrow 0}E\{x^{\varepsilon }(\tau )\}^{2}=2D\tau ,
\end{equation*}
where
\begin{equation}
D=\int_{0}^{\infty }<\xi (\tau +s)\xi (\tau )>ds.
\end{equation}

\renewcommand{\theequation}{B-\arabic{equation}} \setcounter{equation}{0}

\section*{Appendix B}

In this Appendix we derive the equation (\ref{PDE}) for the covariance
\begin{equation}
R\left( \tau ,S,Y\right) =<Z\left( \tau ,S\right) Z\left( \tau ,Y\right) >.
\end{equation}
First, let us find the derivative of $R\left( \tau ,S,Y\right) $ with
respect to time
\begin{equation}
\frac{\partial R}{\partial \tau }=<Z\left( \tau ,S\right) \frac{\partial
Z\left( \tau ,Y\right) }{\partial \tau }+Z\left( \tau ,Y\right) \frac{%
\partial Z\left( \tau ,S\right) }{\partial \tau }>.  \label{deriv}
\end{equation}
Substitution of the derivatives $\frac{\partial Z\left( \tau ,S\right) }{%
\partial \tau }$ and $\frac{\partial Z\left( \tau ,Y\right) }{\partial \tau }
$ from equation (\ref{limitequation}) into (\ref{deriv}) and
averaging give
\begin{eqnarray}
\frac{\partial R}{\partial \tau } &=&\frac{1}{2}\sigma ^{2}S^{2}\frac{%
\partial ^{2}R}{\partial S^{2}}+\frac{1}{2}\sigma ^{2}Y^{2}\frac{\partial
^{2}R}{\partial Y^{2}}+r(S\frac{\partial R}{\partial S}+Y\frac{\partial R}{%
\partial Y}-2R)+  \notag \\
&&\left( S\frac{\partial V_{BS}}{\partial S}-V_{BS}\right) \left\langle
Z\left( \tau ,Y\right) \eta (\tau )\right\rangle +\left( Y\frac{\partial
V_{BS}}{\partial Y}-V_{BS}\right) \left\langle Z\left( \tau ,S\right) \eta
(\tau )\right\rangle .
\end{eqnarray}
This equation involves the correlation functions $\left\langle
Z\left( \tau ,Y\right) \eta (\tau )\right\rangle $ and
$\left\langle Z\left( \tau ,S\right) \eta (\tau )\right\rangle $
that can be found as follows. Since the random process $\eta (\tau
)$ is Gaussian, one can use the Furutsu-Novikov formula to find
\begin{equation}
\left\langle Z\left( \tau ,Y\right) \eta (\tau )\right\rangle
=\int \left\langle \eta (\tau )\eta (\tau _{1})\right\rangle
\times \left\langle \frac{\delta Z\left( \tau ,Y\right) }{\delta
\eta (\tau _{1})}\right\rangle d\tau _{1}  \label{Novikov}
\end{equation}
(see \cite{Moss}). By using delta-correlated function for $\eta
(\tau )$ given by (\ref{delta}), and equation
(\ref{limitequation}), we can find the variational derivative

\begin{equation}
\frac{\delta Z\left( \tau ,Y\right) }{\delta \eta (\tau )}=Y\frac{\partial
V_{BS}}{\partial Y}-V_{BS}\left( \tau ,Y\right) ,
\end{equation}
(see \cite {Fe}) and then the correlation function (\ref{Novikov})
\begin{equation}
\left\langle Z\left( \tau ,Y\right) \eta (\tau )\right\rangle =D\left( Y%
\frac{\partial V_{BS}}{\partial Y}-V_{BS}\left( \tau ,Y\right) \right) .
\end{equation}
The equation for the covariance $R$ becomes
\begin{eqnarray}
\frac{\partial R}{\partial \tau } &=&\frac{1}{2}\sigma ^{2}S^{2}\frac{%
\partial ^{2}R}{\partial S^{2}}+\frac{1}{2}\sigma ^{2}Y^{2}\frac{\partial
^{2}R}{\partial Y^{2}}+r(S\frac{\partial R}{\partial S}+Y\frac{\partial R}{%
\partial Y}-2R)+  \notag \\
&&2D\left( S\frac{\partial V_{BS}}{\partial S}-V_{BS}\left( \tau ,S\right)
\right) \left( Y\frac{\partial V_{BS}}{\partial Y}-V_{BS}\left( \tau
,Y\right) \right) .
\end{eqnarray}

\renewcommand{\theequation}{C-\arabic{equation}} \setcounter{equation}{0}

\section*{Appendix C}

In this Appendix we briefly discuss the derivation of (\ref{U}). Using (\ref
{solution}), one can average $Z^{2}(\tau ,S)$ as follows
\begin{align}
U(\tau ,S)& =\langle Z^{2}(\tau ,S)\rangle =\int_{0}^{\tau }\int_{0}^{\tau
}\int_{0}^{\infty }\int_{0}^{\infty }G(S,S_{1},\tau ,\tau
_{1})G(S,S_{2},\tau ,\tau _{2})  \label{stef} \\
& \left( \frac{\partial V_{BS}}{\partial S_{1}}S_{1}-V_{BS}\left( \tau
_{1},S_{1}\right) \right) \left( \frac{\partial V_{BS}}{\partial S_{2}}%
S_{2}-V_{BS}\left( \tau _{2},S_{2}\right) \right) \langle \eta \left( \tau
_{1}\right) \eta (\tau _{2})\rangle dS_{1}dS_{2}d\tau _{1}d\tau _{2},  \notag
\end{align}
where $\langle \eta \left( \tau _{1}\right) \eta (\tau _{2})\rangle
=2D\delta \left( \tau _{1}-\tau _{2}\right) $, and the intensity of the
noise $D$ is determined by (\ref{dif}). It follows from the property of the
Dirac delta function that
\begin{equation}
U\left( \tau ,S\right) =2D\int_{0}^{\tau }[\int_{0}^{\infty }G(S,S_{1},\tau
,\tau _{1})\left( \frac{\partial V_{BS}}{\partial S_{1}}S_{1}-V_{BS}\left(
\tau _{1},S_{1}\right) \right) dS_{1}]^{2}d\tau _{1}.
\end{equation}

\newpage

\begin{figure}[p]
\centering
\includegraphics[scale=0.5]{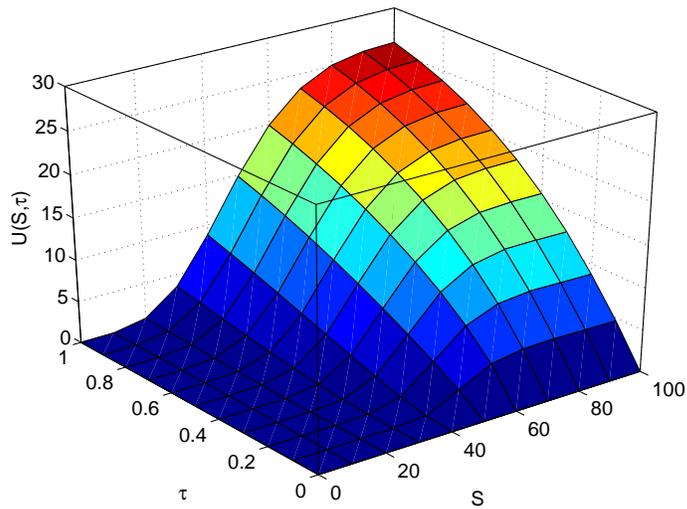} ~
\caption{Variance, $U(\protect\tau ,S)$, with respect to asset
price, $S$, and time $\protect\tau $. The option strike price $K$
is $20$, the volatility $\protect\sigma $ is $0.4$, the interest
rate $r$ is $0.1$, and the constant $D$ is $0.1$.} \label{figure
1}

\end{figure}

\begin{figure}[p]
\centering
\includegraphics[scale=0.8]{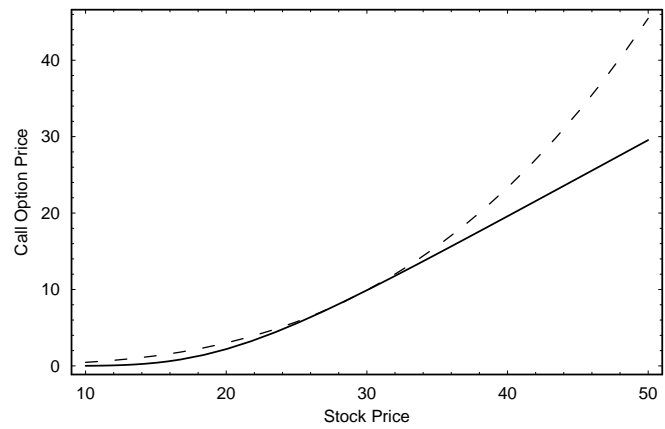} ~
\caption{Effective option price (dashed line) and Black-Scholes
price
(continuous line). Here $\protect\varepsilon =0.1,$ and constants $K$, $%
\protect\sigma $, $r$, and $D$ are taken as in Figs. 1. }
\label{figure 2}
\end{figure}

\begin{figure}[p]
\centering
\includegraphics[scale=0.8]{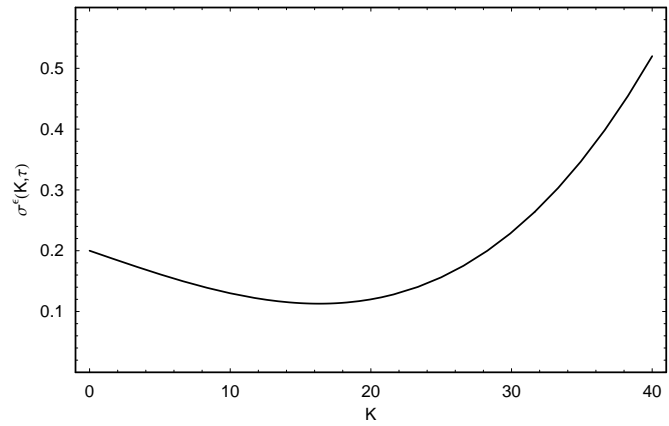} ~
\caption{''Smile'' curve: the implied volatility $ \sigma
^{\varepsilon }\left( K,\tau \right)$ as a function of the strike
price $K$. Here $\protect\varepsilon =0.1$ and constants  $\protect\sigma $, $r$%
, and $D$ are taken as in Fig. 1.} \label{figure 3}
\end{figure}

\end{document}